\documentclass[30pt]{article}  

\usepackage{amsmath}
\usepackage{mathtools}
\usepackage{ntheorem}
\usepackage[margin=1.2in]{geometry}
\usepackage{setspace}
\usepackage{amssymb}
\newtheorem{theorem}{Theorem}[section]
\newtheorem{lemma}{Lemma}[section]

\title{Some coset actions in $G_2(q)$ and distance-transitive graphs} 
\author{Jianxiang Chen} 
\date{\vspace{-5ex}}
\begin{document}

\maketitle

\begin{abstract}
\indent This paper studies whether there are distance-transitive graphs arising from the coset actions of $G_2(q)$ on the subfield subgroup $G_2(\sqrt{q})$ or $G_2(q)$ on the Ree subgroup $^2G_2(q)$. It is found that there are no such graphs, even if the groups are extended by outer automorphisms of $G_2(q)$.
\end{abstract}

\section{Introduction}

\indent \indent
In this paper, a \textit{graph} means a finite connected undirected graph without loops or multiple edges. A graph $\Gamma$ is called \textit{distance-transitive} if for every two pairs $(w,x)$ and $(y,z)$ of the graph such that the distance from $w$ to $x$ is the same as the distance from $y$ to $z$, there is an automorphism of the graph that carries $w$ to $y$ and $x$ to $z$. 

By a result of D.H.Smith$^{[1]}$, distance-transitive graphs can be classified into either \textit{primitive} or non-primitive, and non-primitive distance-transitive graphs can be derived from primitive distance-transitive graphs, so the classification of primitive distance-transitive graphs serves as a step to classify all distance-transitive graphs. 

Praeger, Saxl and Yokoyama$^{[2]}$ showed that a primitive distance-transitive graph is either a Hamming graph, has an automorphism group of affine type, or has an almost-simple automorphism group (a group $G$ is \textit{almost simple} if there is a nonabelian simple group $H$ such that $H \subset G \subset \text{Aut}(H)$; here $H$ is the \textit{socle} of $G$). The primitive distance-transitive graphs having automorphism groups of affine type are classified by John van Bon$^{[3]}$.

The vertex stabilizer in a primitive distance-transitive graph $\Gamma$ is a maximal subgroup of $\text{Aut}\Gamma$. Thanks to the classification of all finite simple groups and their maximal subgroups, the classification of primitive distance-transitive graphs with almost-simple automorphism groups can be dealt with on a case-by-case basis. A lot of work has already been done; for example, the graphs $\Gamma$ where $\text{Aut}\Gamma$ have alternating, sporadic or linear socle are completely determined in [4], [5] and [6] respectively.

There is a combinatorial generalization of distance-transitive graphs called \textit{distance-regular graphs}: a distance-regular graph is a regular graph such that for any two vertices $x$ and $y$, the number of vertices at distance $i$ from $x$ and at distance $j$ from $y$ is independent from $x$ and $y$, and only dependent upon $i$, $j$ and the distance between $x$ and $y$. It is evident from the definition that every distance-transitive graph is distance-regular, because every ordered pair $(x,y)$ of vertices with the same distance between $x$ and $y$ are equivalent. Distance-regular graphs are usually characterized by intersection arrays.

In this paper we investigate two open cases where $\text{Aut}\Gamma$ has socle $G_2(q)$: they correspond to the vertex stabilizer having type $G_2(\sqrt{q})$ and $^2G_2(q)$. More formally, we assume that $\text{Aut}\Gamma=G_2(q)\text{:}X$ where $X$ is a subgroup of the outer automorphism group of $G_2(q)$. The vertex stabilizer of $\Gamma$ is $G_2(\sqrt{q})\text{:}X$ and $^2G_2(q)\text{:}X$, respectively. So the vertices of $\Gamma$ can be identified with the cosets $(G_2(q)\text{:}X)/(G_2(\sqrt{q})\text{:}X)$ and $(G_2(q)\text{:}X)/(^2G_2(q)\text{:}X)$, respectively.

\section{$G_2(q)$ on $G_2(\sqrt{q})$}
\indent \indent Let $G$ be the finite simple group $G_2(q)$, where $q=p^{2n}$, and $H$ the subgroup of $G$ isomorphic to the subfield subgroup $G_2(\sqrt{q})$; all such subgroups are pairwise conjugate in $G$$^{[7,8]}$. Let $X$ be a subgroup of the outer automorphism group of $G$. Let $r=\sqrt{q}$. Let $\sigma$ be the Frobenius isomorphism $\sigma: x \rightarrow x^{r}$. 

 The first thing to note is that $q$ is a power of $3$. If $\Gamma$ is distance-transitive, the permutation character of $\text{Aut}\Gamma$ acting on the vertices of $\Gamma$ is multiplicity-free$^{[3]}$. The multiplicity-free actions of $(G\text{:}X)/(H\text{:}X)$ are determined in [9]: The action is multiplicity-free iff $q$ is a power of $3$ and $X$ contains the graph automorphism of $G_2(q)$. So $q$ is a power of $3$.

 Ross Lawther has computed the suborbit lengths of $G$ acting on $G/H$ when $q$ is a power of $3$$^{[10]}$, shown in Table 1, where $z$ is the $G_2(\bar{\mathbb{F}}_q)$-class representative of $x^{-1}\sigma(x)$ in $G_2(q)$ listed in [11].

The elements of the form $x_{ma+nb}(1)$ are unipotent elements of the algebraic group $G_2(\bar{\mathbb{F}}_q)$, so they have order $3$. Note that $x_{2a+b}(1)$ and $x_{3a+2b}(1)$ commute$^{[10]}$, so $x_{2a+b}(1)x_{3a+2b}(1)$ also has order $3$.

Let $T$ be a Cartan subgroup of the algebraic group $G_2(\bar{\mathbb{F}}_q)$. The Dynkin diagram with respect to $T$ is of type $G_2$. Denote its short simple root as $a$ and long simple root as $b$ under a choice of simple roots. Then $z_1=a+b$, $z_2=a$, $z_3=-2a-b$ are three short roots of the root diagram. The elements $h(x_1,x_2,x_3)$$(x_1,x_2,x_3 \in \bar{\mathbb{F}}_q^*)$ are elements in $T$, mapping $z_i$ to $x_i$, $i=1,2,3$$^{[11]}$.

 The letters $\sigma$, $\tau$, $\theta$, $\eta$ and $\gamma$ in Table 1 represents certain elements of the multiplicative group $\mathbb{F}_{q^3}^*$. Let $\kappa$ be a generator of the multiplicative group $\mathbb{F}_{q^3}^*$, and
\begin{table}[htb]
\centering
\begin{tabular}{lll}
$\sigma = \kappa^{(r+1)(r^3-1)}$, & $\tau=\kappa^{(r-1)(r^3+1)}$, & $\theta=\kappa^{q^2+q+1}$, \\
$\eta=\theta^{r-1}$, & $\gamma=\theta^{r+1}$. &\\
\end{tabular}
\end{table}

The shorthand $h_*(x_1,x_2,x_3)$ represents $h(*^{x_1},*^{x_2},*^{x_3})$ where $*$ is any letter.

\begin{table}[htb]
\centering
\begin{tabular}{ |l|l|l| } 
\hline
$z$ & Size of suborbits &  Number of classes \\
\hline 
$1$ & $1$ & $1$ \\
$x_{3a+2b}(1)$ & $(r^6-1)$ & $1$\\
$x_{2a+b}(1)$ & $(r^6-1)$ & $1$\\
$x_{2a+b}(1)x_{3a+2b}(1)$ & $(r^6-1)(r^2-1)$ & $1$\\
$x_{a+b}(1)x_{3a+b}(1)$ & $r^2(r^6-1)(r^2-1)/2$ & $1$\\
$$ & $r^2(r^6-1)(r^2-1)/2$ & $1$\\
$x_{a}(1)x_{b}(1)$ & $r^4(r^6-1)(r^2-1)$ & $1$\\
& &\\
$h(-1,-1,1)$& $r^4(r^4+r^2+1)$ & 1\\
$h(-1,-1,1)x_b(1)$& $r^4(r^6-1)$ &1\\
$h(-1,-1,1)x_{2a+b}(1)$& $r^4(r^6-1)$ &1\\
$h(-1,-1,1)x_b(1)x_{2a+b}(1)$& $r^4(r^6-1)(r^2-1)/2$ &1\\
$$&$r^4(r^6-1)(r^2-1)/2$ &1\\
& &\\
$h_\gamma(i,-2i,i)$ & $r^5(r^3-1)(r^2-r+1)$ & $(r-3)/2$ \\
$h_\gamma(i,-2i,i)x_{3a+2b}(1)$ & $r^5(r^6-1)(r-1)$ & $(r-3)/2$ \\
$h_\gamma(i,-i,0)$ & $r^5(r^3-1)(r^2-r+1)$ & $(r-3)/2$ \\
$h_\gamma(i,-i,0)x_{2a+b}(1)$ & $r^5(r^6-1)(r-1)$ & $(r-3)/2$ \\
$h_\gamma(i,j,-i-j)$ & $r^6(r^3-1)(r^2-r+1)(r-1)$ & $(r^2-8r+15)/12$ \\
& &\\
$h_\eta(i,-2i,i)$ & $r^5(r^3+1)(r^2+r+1)$ & $(r-1)/2$ \\
$h_\eta(i,-2i,i)x_{3a+2b}(1)$ & $r^5(r^6-1)(r+1)$ & $(r-1)/2$ \\
$h_\eta(i,-i,0)$ & $r^5(r^3+1)(r^2+r+1)$ & $(r-1)/2$ \\
$h_\eta(i,-i,0)x_{2a+b}(1)$ & $r^5(r^6-1)(r+1)$ & $(r-1)/2$ \\
$h_\eta(i,j,-i-j)$ & $r^6(r^3+1)(r^2+r+1)(r+1)$ & $(r^2-4r+3)/12$ \\
& & \\
$h_\theta(i,(r-1)i,-ri)$ & $r^6(r^6-1)$ & $(r-1)^2/4$ \\
$h_\theta(i,ri,-(r+1)i)$ & $r^6(r^6-1)$ & $(r-1)^2/4$ \\
$h_\tau(i,ri,r^2i)$      & $r^6(r^3-1)(r^2-1)(r+1)$ & $r(r+1)/6$ \\ 
$h_\sigma(i,-ri,r^2i)$   & $r^6(r^3+1)(r^2-1)(r-1)$ & $r(r-1)/6$ \\ 

\hline
\end{tabular}
\caption{The suborbits of the coset action $G/H$ when $q$ is a power of $3$.}
\end{table}

\newpage

 The second thing to note is that the subgroup $G_2(\sqrt{q})$ is the fixed subgroup of the Frobenius isomorphism $\sigma: x \rightarrow x^{\sqrt{q}}$. If $\sigma \in X$, then the subgroup $H\text{:}X$ is an involution centralizer in $G\text{:}X$, and the action of $G\text{:}X$ on cosets $(G\text{:}X)/(H\text{:}X)$ can be identified with the action of $G\text{:}X$ on the conjugacy class of $\sigma$ by conjugation. Then it would be possible to apply the methodology of [12], summarized in the following theorem:

\begin{theorem}
Let $\Gamma$ be a distance-transitive graph with distance-transitive group $G$. Suppose that the vertex set $V\Gamma$ of $\Gamma$ is a conjugacy class of involutions in $G$, that $G$ acts on $\Gamma$ by conjugation and that there are elements in $V\Gamma$ which commute in $G$. Take $x,y \in \Gamma$ with $x$ adjacent to $y$. Then at least one of the following statements holds.

· $\Gamma$ is a polygon or an antipodal 2-cover of a complete graph.

· $G$ is a 2-group.

· The order of $xy$ is an odd prime, if $a,b \in \Gamma$ with $ab$ of order $2$, then $a$ and $b$ have maximal distance in $\Gamma$, and if $a,b \in \Gamma$ the order of $ab$ is not 4.

· The elements $x$ and $y$ commute, and if $z \in O_2 (x)$ then $xz$ has order $2$, $4$ or an odd
prime. Moreover either $O_2(C_G(x)) = \langle x\rangle$ or $C_G(x)$ contains a normal subgroup
generated by $p$-transpositions.
\end{theorem}

 By the following lemma of [12], we can assume $\sigma \in X$:

\begin{lemma}

Let $\Gamma$ be a graph on which $G$ acts primitively distance-transitively, and denote by $H$ the stabilizer in $G$ of a vertex of $\Gamma$. Suppose $\sigma$ is an automorphism of $G$.

· If $\sigma$ centralizes $H$ and diam $\Gamma \geq 3$, then $\sigma \in \text{Aut}(\Gamma)$;

· If $\sigma$ normalizes $H$ and diam $\Gamma \geq 5$, then the same conclusion holds. 

\end{lemma}

The suborbit lengths shown in Table 1 indicates that a distance-transitive $\Gamma$ cannot have diameter $2$ (because there are at least $3$ different nontrivial suborbit lengths in $G$ acting on $G/H$, and the outer automorphism group can only fuse together suborbits of the same length). So we may assume $\Gamma$ has diameter $\geq 3$ and thus $\sigma \in X$.

In the rest of this section, we will denote the elements of $G\text{:}X$ in external semidirect product notation; that is, an element of $G\text{:}X$ is written as $(x,y)$, where $x \in G$ and $y \in X$, and the multiplication rule is $(w,x)(y,z)=(wx(y),xz)$. Thus, the conjugation of $(1,\sigma)$ by $(g,1)$ is $(g^{-1},1)(1,\sigma)(g,1)=(g^{-1},\sigma)(g,1)=(g^{-1}\sigma(g),\sigma)$. Also $(g^{-1}\sigma(g),\sigma)(1,\sigma)=(g^{-1}\sigma(g),1)$. If $(1,\sigma)$ commutes with $(g^{-1}\sigma(g),\sigma)$, we have $g^{-1}\sigma(g)=\sigma(g^{-1}\sigma(g))$. And that means $(g^{-1}\sigma(g))^2=g^{-1}\sigma(g)g^{-1}\sigma(g)=g^{-1}\sigma(g)\sigma(g^{-1}\sigma(g))=1$. The reverse implication also holds, in the sense that $(g^{-1}\sigma(g))^2=1$ implies $(1,\sigma)$ commutes with $(g^{-1}\sigma(g),\sigma)$.

There exists elements $g\in G$ such that $g^{-1}\sigma(g)$ is conjugate to $h(-1,-1,1)$. In this case, $g^{-1}\sigma(g)$ has order $2$, so $(1,\sigma)$ commutes with $(g^{-1}\sigma(g),\sigma)$. Thus the assumptions of Theorem 2.1 holds.

To find an element with the form $(g^{-1}\sigma(g),\sigma)$ that is connected to $(1,\sigma)$, the following theorem$^{[3]}$ is applied:

\begin{theorem}
Let $G$ be a primitive distance-transitive group of automorphisms of $\Gamma$ (having diameter $d$) and $x \in V\Gamma$ . Then among the nontrivial $G_x$-orbit lengths, $|\Gamma_1(x)|$ $(\Gamma_i(x)$ means the vertices of $\Gamma$ having distance $i$ to $x)$ is among the two smallest. Moreover, if $|\Gamma_1(x)|$ is not the smallest, then $|\Gamma_d(x)|$ is.
\end{theorem}

Since outer automorphisms can only fuse suborbits with the same length, the smallest suborbits must be one labeled by $x_{3a+2b}(1)$, $x_{2a+b}(1)$ or $x_{2a+b}(1)x_{3a+2b}(1)$. In either case, the element $g^{-1}\sigma(g)$ has order $3$, so $(1,\sigma)$ and $(g^{-1}\sigma(g),\sigma)$ does not commute. So the last case of Theorem 2.1 does not hold. As $G\text{:}X$ is not a 2-group and $\Gamma$ has degree and diameter at least $3$, the first and second case does not hold. It can be concluded that if $\Gamma$ is distance-transitive, then the third case of Theorem 2.1 holds. Specifically, there are no $g\in G\text{:}X$ such that $g^{-1}\sigma(g)$ has order $4$.

There exists $g\in G$ such that $g^{-1}\sigma(g)$ is conjugate to $h_{\gamma}(i, -2i, i)$ when $r \geq 9$, and there exists $g\in G$ such that $g^{-1}\sigma(g)$ is conjugate to $h_{\eta}(i,-2i,i)$. The elements $\gamma$ and $\eta$ has order $r-1$ and $r+1$, respectively. As one of $r-1$ and $r+1$ is divisible by $4$, it is possible to choose $i$ such that one of $\gamma^i$ and $\eta^i$ has order $4$ in the multiplicative group $\mathbb{F}_{q^3}^*$. But this means that $g^{-1}\sigma(g)$ has order $4$.

So there is no primitive distance-transitive graph with automorphism group $G_2(q)\text{:}X$ and vertex stabilizer $G_2(\sqrt{q})\text{:}X$.

\section{$G_2(q)$ on $^2G_2(q)$}
\indent \indent Let $G$ be the finite simple group $G_2(q)$, where $q=3^{2n+1}$ ($n\geq 0$), and $H$ the subgroup of $G$ isomorphic to the Ree group $^2G_2(q)$; all such subgroups are pairwise conjugate in $G$$^{[7]}$. Let $X$ be a subgroup of the outer automorphism group of $G$.

 Ross Lawther has computed the suborbit lengths of $G$ acting on $G/H$$^{[10]}$:

\begin{table}[htb]
\centering
\begin{tabular}{ |l|l| } 
\hline
Size of suborbits &  Number of classes \\
\hline 
$1$ & $1$ \\
$(q^3+1)(q-1)$ & $1$\\
$q(q^3+1)(q-1)/2$ & $1$\\
$q(q^3+1)(q-1)/2$ & $1$\\
$q^2(q^3+1)(q-1)$ & $1$\\
&\\
$q^2(q^2-q+1)$ & $1$\\
$q^2(q^3+1)(q-1)/2$ & $1$\\
$q^2(q^3+1)(q-1)/2$ & $1$\\
&\\
$q^3(q^3+1)$ & $(q-3)/2$ \\
$q^3(q^2-q+1)(q-1)$ & $(q-3)/6$ \\
$q^3(q^2-1)(q-3m+1)$ & $(q-3m)/6$ \\
$q^3(q^2-1)(q+3m+1)$ & $(q+3m)/6$ \\
\hline
\end{tabular}
\caption{The suborbits of the coset action $G/H$. Note: $q=3m^2$}
\end{table}

Note that three rows of this table would be absent if $q=3$; but in this case, the graph would have $|G|/|H|=4245696/1512=2808$ vertices and diameter at least $6$ (because the outer automorphism group can only fuse together suborbits of the same length). A graph of this size would be covered in Chapter 14. of [13], but [13] does not contain any intersection arrays of a distance-regular graph with $2808$ vertices and diameter at least $6$. So we will assume $q \geq 27$ in the rest of this section.

An inequality from [14] states that a distance-regular graph of diameter $d$ and $v$ vertices satisfies $d < 8/3 \log_2(v)$. If there were a distance-regular graph arising from the group action, it would have $q^3(q^3-1)(q+1)$ vertices and diameter at least $(q+6)/|X|=(q+6)/2(2n+1)$, because the outer automorphism group of $G$ is a direct product of the graph automorphism (of order $2$) and the field automorphism (of order $2n+1$), and it can never fuse more than $|X|$ suborbits into a suborbit. This inequality is only satisfied when $n \leq 3$.

The main theorem to deal with the cases when $n \leq 3$ is the following theorem from [12]:
\begin{theorem}
Let $\Gamma$ be a graph of diameter $d$ on which the group $G$ acts distance-transitively as a group of automorphisms. For a vertex $x \in \Gamma$, denote by $G_x^i$ the kernel of the action of the stabilizer of $x$ in $G$ on $\Gamma_i(x)$. If, for some $i\geq1$, we have $G_x^i \neq 1$, then $G^i_x \subsetneq G^{i-1}_x  \subsetneq ...  \subsetneq G^1_x$ or $G^i_x \subsetneq G^{i+1}_x  \subsetneq ...  \subsetneq G^{d}_x$.
\end{theorem}

Since all the suborbit lengths are proper divisors of $|H|=q^3(q^3+1)(q-1)$, the kernel of $H$ acting on any of the suborbits is nontrivial, and so does the extension $H\text{:}X$.

The last two rows of Table 2 contain two suborbits. The respective kernel size (the size of $G_x^i$) of one row is divisible by $19$ when $n=1$, $31$ when $n=2$, $43$ when $n=3$, and the other divisible by $37$ when $n=1$, $271$ when $n=2$, $2269$ when $n=3$. By Theorem 2.2, the suborbit corresponding to $\Gamma_1(x)$ has length $(q^3+1)(q-1)$ or $q^2(q^2-q+1)$. Thus $G^1_x$ is not divisible by any of the two kernel sizes, and it is impossible to have $G^i_x \subsetneq G^1_x$ for any of the two suborbits. So $G^i_x \subsetneq G^d_x$ for both suborbits, and $G^d_x$ would have order divisible by the primes $19\times37$ when $n=1$, $31\times271$ when $n=2$, and $43\times2269$ when $n=3$. But no suborbit has such a kernel, even if $X$ is nontrivial (outer automorphisms can only multiply the order of the kernels by some number containing the prime factors $2$, $3$, $5$ and $7$.).

So there is no primitive distance-transitive graph with automorphism group $G_2(q)\text{:}X$ and vertex stabilizer $^2G_2(q)\text{:}X$.

\section*{References}

\indent \indent 
1. Smith, Derek H. "Primitive and imprimitive graphs." The Quarterly Journal of Mathematics 22.4 (1971): 551-557.

2. Praeger, Cheryl E., Jan Saxl, and Kazuhiro Yokoyama. "Distance transitive graphs and finite simple groups." Proceedings of the London Mathematical Society 3.1 (1987): 1-21.

3. van Bon, John. "Finite primitive distance-transitive graphs." European Journal of Combinatorics 28.2 (2007): 517-532.

4. Liebeck, Martin W., Cheryl E. Praeger, and Jan Saxl. "Distance transitive graphs with symmetric or alternating automorphism group." Bulletin of the Australian Mathematical Society 35.1 (1987): 1-25.

5. Ivanov, A. A., et al. "Distance-transitive representations of the sporadic groups." Communications in Algebra 23.9 (1995): 3379-3427.

6. van Bon, John, and Arjeh M. Cohen. "Linear groups and distance-transitive graphs." European Journal of Combinatorics 10.5 (1989): 399-411.

7. Kleidman, Peter B. "The maximal subgroups of the Chevalley groups $G_2(q)$ with $q$ odd, the Ree groups $^2G_2(q)$, and their automorphism groups." Journal of Algebra 117.1 (1988): 30-71.

8. Cooperstein, Bruce N. "Maximal subgroups of $G_2(2^n)$." Journal of Algebra 70.1 (1981): 23-36.

9. Lawther, R. "Some (almost) multiplicity-free coset actions." Groups, Combinatorics and Geometries, in: LMS Lecture Note Series 165 (1992).

10. Lawther, R. "Some coset actions in $G_2(q)$." Proceedings of the London Mathematical Society 3.1 (1990): 1-17.

11. Enomoto, Hikoe. "The conjugacy classes of Chevalley groups of type $G_2$ over finite fields of characteristic 2 or 3." Journal of the Faculty of Science, University of Tokyo. Sect. 1, Mathematics, astronomy, physics, chemistry 16.3 (1970): 497-512.

12. van Bon, John. "On distance-transitive graphs and involutions." Graphs and Combinatorics 7.4 (1991): 377-394.

13. A. E. Brouwer, A. M. Cohen, and A. Neumaier, Distance-Regular Graphs, Ergebnisse der Mathematik 3.18, Springer-Verlag, Heidelberg, (1989).

14. Bang, Sejeong, Akira Hiraki, and Jacobus H. Koolen. "Improving diameter bounds for distance-regular graphs." European Journal of Combinatorics 27.1 (2006): 79-89.

\end{document}